\documentclass[a4paper,10pt]{article}
\usepackage{amsmath}
\usepackage{amsfonts}
\usepackage{amssymb}
\usepackage{amscd}
\usepackage{pb-diagram}
\usepackage{epic}
\usepackage[french]{babel}
\usepackage{url}
\usepackage{color}
\title{La Grassmannienne non-lin\'eaire comme vari\'et\'e fr\'ech\'etique homog\`ene}
\author{}
\date{}
\begin{document}

\newtheorem{lemme}{Lemme}[section]
\newtheorem{definition}[lemme]{D\'efinition}
\newtheorem{proposition}[lemme]{Proposition}
\newtheorem{corollaire}[lemme]{Corollaire}
\newtheorem{theoreme}[lemme]{Th\'eor\`eme}
\newtheorem{remarque}[lemme]{Remarque}
\bibliographystyle{alpha}
\newenvironment{maliste}%
{ \begin{list}%
	{$\bullet$}%
	{\setlength{\labelwidth}{30pt}%	 \setlength{\leftmargin}{35pt}%
	 \setlength{\itemsep}{\parsep}}}%
{ \end{list} }
\maketitle

\begin{center}
\begin{large}Molitor Mathieu\end{large}\\
\textit{Facult\'e des sciences de base -- Institut de math\'ematiques B}\\
\textit{Ecole Polytechnique F\'ed\'erale de Lausanne}\\
\textit{1015 Lausanne, Switzerland} \\  
\textit{e-mail:}\,\,\url{pergame.mathieu@gmail.com}
\end{center}
\begin{abstract}
Soit $(M,g)$ une vari\'et\'e riemannienne compacte de dimension $n\,.$ Pour $k\in \{0,...,n\}$, notons $Gr_{k}(M)$ l'ensemble des sous-vari\'et\'es compactes, connexes, orient\'ees de $M$ et de dimension $k$. Cet ensemble est appel\'e la Grassmannienne non-lin\'eaire. Dans cet article, nous munissons $Gr_{k}(M)$ d'une structure de vari\'et\'e fr\'ech\'etique et d\'eveloppons les propri\'et\'es les plus imm\'ediates de cette vari\'et\'e. Notamment, si $\Sigma\in Gr_{k}(M)$, nous montrons que $\text{Emb}(\Sigma,M)$, l'espace des plongements de $\Sigma$ dans $M$, est l'espace total d'un fibr\'e principal ayant pour base la r\'eunion de certaines composantes connexes de $Gr_{k}(M)\,.$ Nous montrons aussi que les composantes connexes de $Gr_{k}(M)$ sont homog\`enes sous l'action naturelle du groupe des diff\'eomorphismes de $M\,.$
\end{abstract}

\section*{Introduction}
Pour une vari\'et\'e $M$ donn\'ee, l'\'etude des sous-vari\'et\'es de $M\,,$ du point de vue de la g\'eom\'etrie de dimension infinie, am\`ene tr\`es naturellement \`a consid\'erer l'ensemble des sous-vari\'et\'es compactes, connexes et orient\'ees de $M$ comme une vari\'et\'e fr\'ech\'etique. Cet ensemble est appel\'e la Grassmannienne non-lin\'eaire et est not\'ee $Gr_{k}(M)$ ($k$ \'etant la dimension des sous-vari\'et\'es consid\'er\'ees)\,. Il y a essentiellement deux fa\c{c}ons d'appr\'ehender $Gr_{k}(M)\,:$
\begin{center}
\begin{maliste}
\item la premi\`ere, qui est amplement d\'evelopp\'ee dans \cite{Kriegl-Michor}, consiste \`a prendre $\Sigma\in Gr_{k}(M)$ et \`a consid\'erer l'espace des plongements $\text{Emb}\,(\Sigma,M)\,.$ On peut alors identifier $Gr_{k}(M)$ (ou plut\^ot la r\'eunion de certaines composantes connexes de $Gr_{k}(M)$) comme \'etant le quotient de $\text{Emb}\,(\Sigma,M)$ par rapport \`a l'action naturelle du groupe $\text{Diff}^{+}(\Sigma)$  des diff\'eomorphismes de $\Sigma$ qui pr\'eservent une forme de volume donn\'ee, sur $\text{Emb}\,(\Sigma,M)\,.$ Par construction m\^eme, on obtient ainsi une structure lisse de fibr\'e principal $\text{Diff}^{+}(\Sigma)\hookrightarrow \text{Emb}\,(\Sigma,M)\rightarrow \text{Emb}\,(\Sigma,M)/\text{Diff}^{+}(\Sigma)$ (voir Theorem 44.1., page 474 de \cite{Kriegl-Michor}).\\
\item La deuxi\`eme approche, plus intuitive, consiste \`a modeler directement $Gr_{k}(M)$ sur des espaces de sections $\Gamma_{C^{\infty}}(\Sigma,N\Sigma)$ o\`u $\Sigma\in Gr_{k}(M)$ et $N\Sigma$ d\'esigne le fibr\'e normal de $\Sigma$ dans $M\,.$  Cette approche est esquiss\'ee dans \cite{Hamilton}\,.
\end{maliste}
\end{center}
La premi\`ere partie de cet article s'attache \`a d\'ecrire tr\`es explicitement la construction \'ebauch\'ee par Hamilton  dans la cat\'egorie des vari\'et\'es fr\'ech\'etiques mod\'er\'ees (``tame'' en anglais, voir \cite{Hamilton}). La deuxi\`eme partie fait le lien entre les deux points de vue cit\'es. Nous y montrons notamment un th\'eor\`eme analogue au Theorem 44.1. de \cite{Kriegl-Michor}\,, c'est-\`a-dire, nous montrons que l'espace des plongements $\text{Emb}\,(\Sigma,M)$ est un fibr\'e principal, de groupe de structure $\text{Diff}^{+}(\Sigma)$ et dont la base est une r\'eunion de certaines composantes connexes de $Gr_{k}(M)\,.$ Enfin dans la troisi\`eme partie, nous montrons que les composantes connexes de $Gr_{k}(M)$ sont homog\`enes sous l'action naturelle (et lisse) de $\text{Diff}^{\,0}(M),$ la composante connexe en l'\'el\'ement neutre du groupe des diff\'eomorphismes de $M\,.$ Cette homog\'en\'eit\'e est d\'ej\`a mentionn\'ee, mais admise sans preuve, dans \cite{Ismagilov}\,. En revanche, en adoptant le point de vue de \cite{Kriegl-Michor} qui d\'efinit la Grassmannienne non-lin\'eaire comme le quotient de $\text{Emb}\,(\Sigma,M)$ par l'action de $\text{Diff}^{+}(\Sigma)\,,$ l'homog\'en\'eit\'e des composantes connexes de $\text{Emb}\,(\Sigma,M)$ sous $\text{Diff}^{\,0}(M)$ (homog\'en\'eit\'e qui est une cons\'equence directe d'un r\'esultat classique de topologie diff\'erentielle sur les extensions des isotopies en diff\'eotopies, voir \cite{Hirsch}, theorem 1.3. page 180) implique automatiquement l'homog\'en\'eit\'e des composantes connexes correspondantes du quotient $\text{Emb}\,(\Sigma,M)/\text{Diff}^{+}(\Sigma)\,.$ C'est cette approche qui est utilis\'ee dans \cite{Haller-Vizman}\,. Ici encore, et contrairement \`a cette derni\`ere, nous utilisons l'approche de \cite{Hamilton} et regardons $Gr_{k}(M)$ comme une collection de sous-vari\'et\'es dont la structure diff\'erentielle est celle expliqu\'ee dans la section \ref{premiere partie}. Ce faisant, un travail suppl\'ementaire est n\'ecessaire pour montrer l'homog\'en\'eit\'e des composantes de $Gr_{k}(M)$\,.\\\\
La notion de calcul diff\'erentiel sur un espace de Fr\'echet n'\'etant pas ``canonique'', nous joignons un tr\`es court appendice traitant des deux notions de calcul diff\'erentiel les plus courantes sur un espace de Fr\'echet (celle d\'evelopp\'ee par exemple dans \cite{Hamilton}, et celle utilisant la notion de courbes lisses qui est  d\'evelopp\'ee dans \cite{Kriegl-Michor}). Ces deux notions \'etant identiques (sur un espace de Fr\'echet), nous utiliserons indiff\'eremment l'une ou l'autre dans ce texte.

\section{La structure de vari\'et\'e de $Gr_{k}(M)$}\label{premiere partie}
Pour munir $Gr_{k}(M)$ d'une structure de vari\'et\'e, nous allons construire explicitement un atlas sur $Gr_{k}(M)\,.$ Pour ce faire, prenons $\Sigma\in Gr_{k}(M)$ et introduisons les notations et objets suivants :
\begin{description}
\item[$\bullet$]  $\Theta_{\Sigma}$ un ouvert contenant la section nulle du fibr\'e normal $N\Sigma$ de $\Sigma$ dans $M\,,$ convexe fibre par fibre et tel que l'application
$$
\tau_{\Sigma}\,:\,\Theta_{\Sigma}\rightarrow M\,,v\in N\Sigma_{x}\mapsto\text{exp}_{x}(v)
$$
soit un diff\'eomorphisme de $\Theta_{\Sigma}$ sur son image;\\
\item[$\bullet$]  $\mathcal{U}_{\Sigma}:=\{s\in \Gamma_{C^{\infty}}(\Sigma, N\Sigma)\,\vert\,s(\Sigma)\subseteq\Theta_{\Sigma}\}\,;$\\
\item[$\bullet$]  $\varphi_{\Sigma}\,:\,\mathcal{U}_{\Sigma}\rightarrow Gr_{k}(M)\,,$ application d\'efinie par $\varphi_{\Sigma}(s)= \tau_{\Sigma}\big(s(\Sigma)\big)\,,$ cette derni\`ere sous-vari\'et\'e \'etant munie de l'orientation induite par le diff\'eomor-
\\phisme $\Sigma\rightarrow \tau_{\Sigma}\big(s(\Sigma)\big), x\mapsto \tau_{\Sigma}\big(s(x)\big)\,.$
\end{description}
Montrons que $\{(\varphi_{\Sigma}(\mathcal{U}_{\Sigma}),\varphi_{\Sigma}^{-1})\,\vert\,\Sigma\in Gr_{k}(M)\}$ est un atlas diff\'erentiable de $Gr_{k}(M)$ au moyen des deux lemmes suivants.
\begin{lemme}\label{atlas}
Pour $\Sigma_{1},\Sigma_{2}\in Gr_{k}(M)$, l'ensemble $\varphi_{\Sigma_{1}}^{-1}\big(\varphi_{\Sigma_{1}}(\mathcal{U}_{\Sigma_{1}})\cap\varphi_{\Sigma_{2}}(\mathcal{U}_{\Sigma_{2}})\big)$ est un ouvert de $\Gamma_{C^{\infty}}(\Sigma_{1}, N\Sigma_{1})\,.$
\end{lemme}
\textbf{D\'emonstration.}  Montrons le par l'absurde en supposant que $\varphi_{\Sigma_{1}}^{-1}\big(\varphi_{\Sigma_{1}}(\mathcal{U}_{\Sigma_{1}})\cap\varphi_{\Sigma_{2}}(\mathcal{U}_{\Sigma_{2}})\big)$ 
ne soit pas un ouvert de l'espace m\'etrique $\Gamma_{C^{\infty}}(\Sigma_{1}, N\Sigma_{1})\,.$ On peut alors trouver une section $s\in \varphi_{\Sigma_{1}}^{-1}\big(\varphi_{\Sigma_{1}}(\mathcal{U}_{\Sigma_{1}})\cap\varphi_{\Sigma_{2}}(\mathcal{U}_{\Sigma_{2}})\big)$ et une suite de sections $(s_{n})_{n\in \mathbb{N}}$ telles que $s_{n}\rightarrow s$ dans $\Gamma_{C^{\infty}}(\Sigma_{1}, N\Sigma_{1})$ et $s_{n}\not\in \varphi_{\Sigma_{1}}^{-1}\big(\varphi_{\Sigma_{1}}(\mathcal{U}_{\Sigma_{1}})\cap\varphi_{\Sigma_{2}}(\mathcal{U}_{\Sigma_{2}})\big)$ pour tout $n\in \mathbb{N}\,.$ \\Prenons aussi $U$ un voisinage ouvert de $\Sigma:=\varphi_{\Sigma_{1}}(s)$ inclu dans $\tau_{\Sigma_{1}}(\Theta_{\Sigma_{1}})\cap\tau_{\Sigma_{2}}(\Theta_{\Sigma_{2}})\,.$
L'ouvert $U$ peut \^etre vu simultan\'ement comme une fibration (non-lin\'eaire) au-dessus de $\Sigma_{1}$ et de $\Sigma_{2}$ munie des projections $\pi_{1}$ et $\pi_{2}\,:$
$$
\pi_{i}\,:\,U\rightarrow\Sigma_{i}\,.
$$
Remarquons que $\pi_{i}\vert_{\Sigma}\::\,\Sigma\rightarrow\Sigma_{i}$ est un diff\'eomorphisme et que pour $x\in U$ on a :
$$
\pi_{i}(x)=\pi_{\text{N}\Sigma_{i}}\Big(\tau_{\Sigma_{i}}^{-1}(x)\Big)\,,
$$
o\`u $\pi_{\text{N}\Sigma_{i}}\,:\,N\Sigma_{i}\rightarrow\Sigma_{i}$ est la projection canonique. Nous allons montrer que l'application
$$
\Sigma_{1}\rightarrow\Sigma_{2},\,m\mapsto(\pi_{2}\circ\tau_{\Sigma_{1}}\circ s_{n})(m)
$$
est un diff\'eomorphisme pour $n$ assez grand. 

Remarquons d\'ej\`a que puisque $s_{n}\rightarrow s$ dans $\Gamma_{C^{\infty}}(\Sigma_{1}, N\Sigma_{1}),\,\,\tau_{\Sigma_{1}}\big(s_{n}(\Sigma_{1})\big)\subseteq U$ pour $n$ assez grand, et donc l'application ci-dessus a un sens. Nous allons travailler localement, prenons $x\in \Sigma_{1}\,.$ Puisque $s_{n}\rightarrow s$ dans $\Gamma_{C^{\infty}}(\Sigma_{1}, N\Sigma_{1})$, il existe $y\in \Sigma_{2}$ tel que $(\pi_{2}\circ\tau_{\Sigma_{1}}\circ s_{n})(x)\rightarrow y\,.$ Prenons alors des cartes trivialisantes :
\begin{center}
$
\begin{diagram}
\node{\pi_{N\Sigma_{1}}^{-1}(W)} \arrow[2]{e,t}{\displaystyle\Psi_{W}} \arrow{se,b}{\displaystyle\pi_{N\Sigma_{1}}} \node[2]{W\times\mathbb{R}^{n-k}} \arrow{sw,b}{\displaystyle pr_{1}}\\
\node[2]{W}
\end{diagram}
$
\:\:\:\:
$
\begin{diagram}
\node{\pi_{N\Sigma_{2}}^{-1}(\Omega)} \arrow[2]{e,t}{\displaystyle\Psi_{\Omega}} \arrow{se,b}{\displaystyle\pi_{N\Sigma_{2}}} \node[2]{\Omega\times\mathbb{R}^{n-k}} \arrow{sw,b}{\displaystyle pr_{1}}\\
\node[2]{\Omega}
\end{diagram}
$
\end{center}
avec $x\in W\subseteq\Sigma_{1}\,,y\in \Omega\subseteq\Sigma_{2}$ et telle que $(\pi_{2}\circ\tau_{\Sigma_{1}}\circ s_{n})(W)\subseteq \Omega$ \`a partir d'un certain rang.\\ Nous avons alors :
$$
(\pi_{2}\circ\tau_{\Sigma_{1}}\circ s_{n})(x)=\Big(\underbrace{\pi_{N\Sigma_{2}}\circ\Psi_{\Omega}^{-1}}_{(z,w)\mapsto z}\circ\underbrace{\Psi_{\Omega}\circ\tau_{\Sigma_{2}}^{-1}\circ\tau_{\Sigma_{1}}\circ\Psi^{-1}_{W}}_{(y,v)\mapsto(\tau^{1}(y,v),\tau^{2}(y,v))}\circ\underbrace{\Psi_{W}\circ s_{n}}_{x\mapsto(x,\tilde{s}_{n}(x))}\Big)(x)\,.
$$
Cette application a pour diff\'erentielle :
\begin{eqnarray*}
\left(
Id\,\,0
\right)
\left(
\begin{array}{cc}
 \dfrac{\partial\,\tau^{1}}{\partial\,x}& \dfrac{\partial\,\tau^{1}}{\partial\,y} \\
\dfrac{\partial\,\tau^{2}}{\partial\,x} & \dfrac{\partial\,\tau^{2}}{\partial\,y}
\end{array}
\right)
\left(
\begin{array}{c}
 Id\\
(\tilde{s}_{n})_{*_{x}}
\end{array}
\right)&=&
\left(
Id\,\,0
\right)
\left(
\begin{array}{cc}
 \dfrac{\partial\,\tau^{1}}{\partial\,x}+ \dfrac{\partial\,\tau^{1}}{\partial\,y}(\tilde{s}_{n})_{*_{x}} \\
\dfrac{\partial\,\tau^{2}}{\partial\,x} + \dfrac{\partial\,\tau^{2}}{\partial\,y}(\tilde{s}_{n})_{*_{x}}
\end{array}
\right)\\
&=&\dfrac{\partial\,\tau^{1}}{\partial\,x}+ \dfrac{\partial\,\tau^{1}}{\partial\,y}(\tilde{s}_{n})_{*_{x}}\\
&\rightarrow&\dfrac{\partial\,\tau^{1}}{\partial\,x}+ \dfrac{\partial\,\tau^{1}}{\partial\,y}(\tilde{s})_{*_{x}}\,.
\end{eqnarray*}
La fl\`eche ci-dessus signifie uniquement que nous avons convergence dans un espace de matrices vers la matrice $\frac{\partial\,\tau^{1}}{\partial\,x}+ \frac{\partial\,\tau^{1}}{\partial\,y}(\tilde{s})_{*_{x}}$ qui est inversible puisque cette matrice repr\'esente la diff\'erentielle du diff\'eomorphisme $\big(\pi_{2}\vert_{\Sigma}\big)\circ\big(\pi_{1}\vert_{\Sigma}\big)^{-1}\,:\,\Sigma_{1}\rightarrow\Sigma_{2}\,.$
On en d\'eduit qu'\`a partir d'un certain rang, l'application $\Sigma_{1}\rightarrow\Sigma_{2},\,m\mapsto(\pi_{2}\circ\tau_{\Sigma_{1}}\circ s_{n})(m)$ est partout un diff\'eomorphisme local. Pour montrer que c'est un diff\'eomorphisme globale, il suffit de montrer que cette application est injective. Si ce n'\'etait jamais le cas, pour tout $n\in \mathbb{N}$, on pourrait trouver $x_{n},y_{n}\in \Sigma_{1},\, x_{n}\neq y_{n}$ tels que :
$$
(\pi_{2}\circ\tau_{\Sigma_{1}}\circ s_{n})(x_{n})=(\pi_{2}\circ\tau_{\Sigma_{1}}\circ s_{n})(y_{n})
$$
pour tout $n\in \mathbb{N}.$ Par compacit\'e, nous pouvons supposer que $x_{n}\rightarrow x\in \Sigma_{1}$ et $y_{n}\rightarrow y\in \Sigma_{1}\,.$
En utilisant les semi-normes qui d\'efinissent la topologie de $\Gamma_{C^{\infty}}(\Sigma_{1}, N\Sigma_{1})$ (voir par exemple \cite{dieudonne}, page 236), on constate facilement que :
$$
s_{n}(x_{n})\rightarrow s(x)\,\,\,\,\,\text{et}\,\,\,\,\,s_{n}(y_{n})\rightarrow s(y)
$$
et donc
\begin{eqnarray*}
(\pi_{2}\circ\tau_{\Sigma_{1}}\circ s_{n})(x_{n})&=&(\pi_{2}\circ\tau_{\Sigma_{1}}\circ s_{n})(y_{n})\\
\downarrow\:\:\:\:\:\:\:\:\:\:\:\:\:\:&&\:\:\:\:\:\:\:\:\:\:\:\:\:\:\downarrow\\
\big(\pi_{2}\circ\tau_{\Sigma_{1}}\big)\big(s(x)\big)&=&\big(\pi_{2}\circ\tau_{\Sigma_{1}}\big)\big(s(y)\big)\\
\Rightarrow\:\:\:\:\big(\pi_{2}\vert_{\Sigma}\circ\tau_{\Sigma_{1}}\big)\big(s(x)\big)&=&\big(\pi_{2}\vert_{\Sigma}\circ\tau_{\Sigma_{1}}\big)\big(s(y)\big)\\
\Rightarrow\:\:\:\:\:\:\:\:\:\:\:\:\:\:\:\:\:\:\:\:\:\:\:\:\:\:\:\:\:\:\:\:s(x)&=&s(y)\\
\Rightarrow\:\:\:\:\:\:\:\:\:\:\:\:\:\:\:\:\:\:\:\:\:\:\:\:\:\:\:\:\:\:\:\:\:\:\:\:\:\:x&=&y\,.
\end{eqnarray*}
Ici on a utilis\'e le fait que $\pi_{2}\vert_{\Sigma}\::\,\Sigma\rightarrow\Sigma_{2}$ et $\tau_{\Sigma_{1}}$ sont des diff\'eomorphismes.\\\\
De plus, nous pouvons supposer (voir Appendice, Proposition \ref{special curve}), qu'il existe une courbe lisse $\sigma\,:\,\mathbb{R}\rightarrow\Theta_{\Sigma_{1}}\subseteq\Gamma_{C^{\infty}}(\Sigma_{1}, N\Sigma_{1})$ de $\Gamma_{C^{\infty}}(\Sigma_{1}, N\Sigma_{1})$ telle que $\sigma_{\frac{1}{n}}=s_{n}$ et $\sigma_{0}=s\,.$\\
Consid\'erons alors l'application suivante :
$$
\Lambda\,:\,\mathbb{R}\times\Sigma_{1}\rightarrow\mathbb{R}\times\Sigma_{2}\,,\,(t,x)\mapsto\big(t,\,(\pi_{2}\circ\tau_{\Sigma_{1}}\circ \sigma_{t})(x)\big)\,.
$$
On a que
$$
\Lambda_{*_{(0,x)}}=
\left(
\begin{array}{cc}
\text{Id} & 0 \\
\ast & (\pi_{2}\circ\tau_{\Sigma_{1}}\circ s)_{*_{x}}
\end{array}
\right)
$$
est un isomorphisme puisque l'application $\pi_{2}\circ\tau_{\Sigma_{1}}\circ s=\big(\pi_{2}\vert_{\Sigma}\big)\circ\
\big(\pi_{1}\vert_{\Sigma}\big)^{-1}$ est un diff\'eomorphisme. On en d\'eduit que $\Lambda$ est un diff\'eomorphisme local en $(0,x)\,.$ Mais alors, de part l'\'equivalence suivante :
$$
(\pi_{2}\circ\tau_{\Sigma_{1}}\circ s_{n})(x_{n})=(\pi_{2}\circ\tau_{\Sigma_{1}}\circ s_{n})(y_{n})\,\,\,\,\,\,\Leftrightarrow\,\,\,\,\,\,\Lambda\bigg(\dfrac{1}{n},x_{n}\bigg)=\Lambda\bigg(\dfrac{1}{n},y_{n}\bigg)\,,
$$
il en r\'esulte que pour $n$ assez grand, $x_{n}=y_{n}$, $\Lambda$ devenant injective au voisinage de $(0,x)=(0,y)$, d'o\`u une contradiction.
On en d\'eduit donc que pour $n$ assez grand, $\pi_{2}\circ\tau_{\Sigma_{1}}\circ s_{n}\::\:\Sigma_{1}\rightarrow\Sigma_{2}$ est un diff\'eomorphisme.\\
Ce dernier r\'esultat entraine que $\varphi_{\Sigma_{1}}^{-1}\big(\varphi_{\Sigma_{1}}(\mathcal{U}_{\Sigma_{1}})\cap\varphi_{\Sigma_{2}}(\mathcal{U}_{\Sigma_{2}})\big)$ est un ouvert de $\Gamma_{C^{\infty}}(\Sigma_{1}, N\Sigma_{1})$ car
\begin{eqnarray*}
&&\varphi_{\Sigma_{1}}(s_{n})=\varphi_{\Sigma_{2}}\big(\tau_{\Sigma_{2}}^{-1}\circ\tau_{\Sigma_{1}}\circ(\pi_{2}\circ\tau_{\Sigma_{1}}\circ s_{n})^{-1}\big)\in\varphi_{\Sigma_{2}}\big(U_{\Sigma_{2}}\big)\\
&\Rightarrow&s_{n}\in \varphi_{\Sigma_{1}}^{-1}\big(\varphi_{\Sigma_{1}}(\mathcal{U}_{\Sigma_{1}})\cap\varphi_{\Sigma_{2}}(\mathcal{U}_{\Sigma_{2}})\big)
\end{eqnarray*}
ce qui est une contradiction avec notre hypoth\`ese.$\hfill\square$\\

\begin{lemme}
L'application
$$
\varphi_{\Sigma_{2}}^{-1}\circ\varphi_{\Sigma_{2}}\::\:\varphi_{\Sigma_{1}}^{-1}\big(\varphi_{\Sigma_{1}}(\mathcal{U}_{\Sigma_{1}})\cap\varphi_{\Sigma_{2}}(\mathcal{U}_{\Sigma_{2}})\big)\rightarrow\varphi_{\Sigma_{2}}^{-1}\big(\varphi_{\Sigma_{1}}(\mathcal{U}_{\Sigma_{1}})\cap\varphi_{\Sigma_{2}}(\mathcal{U}_{\Sigma_{2}})\big)
$$
 est lisse mod\'er\'ee (``tame'' en anglais, voir par exemple \cite{Hamilton}).
\end{lemme}
\textbf{D\'emonstration.}  Prenons $\Sigma_{1},\Sigma_{2},\Sigma:=\varphi_{\Sigma_{1}}(s)$ et $U$ comme dans le Lemme \ref{atlas}\,. Nous allons montrer que l'application ci-dessus est lisse mod\'er\'ee sur un voisinage de $s$ dans $\Gamma_{C^{\infty}}(\Sigma_{1}, N\Sigma_{1})\,.$ En fait, nous avons d\'ej\`a vu qu'il existe un voisinage $\mathcal{W}$ de $s$ dans $\Gamma_{C^{\infty}}(\Sigma_{1}, N\Sigma_{1})$ tel que l'application $\sigma\in \mathcal{W} \mapsto \pi_{2}\circ\tau_{\Sigma_{1}}\circ \sigma\in C^{\infty}(\Sigma_{1},\Sigma_{2})$ soit \`a valeurs dans $\text{Diff}(\Sigma_{1},\Sigma_{2})$ et forme ainsi une application lisse mod\'er\'ee de $\mathcal{W}\subseteq \Gamma_{C^{\infty}}(\Sigma_{1}, N\Sigma_{1})$ dans $\text{Diff}(\Sigma_{1},\Sigma_{2})\,.$ Il en r\'esulte que l'application $\mathcal{W}\rightarrow \Gamma_{C^{\infty}}(\Sigma_{2}, N\Sigma_{2}),\,\sigma\mapsto \tau_{\Sigma_{2}}^{-1}\circ\tau_{\Sigma_{1}}\circ(\pi_{2}\circ\tau_{\Sigma_{1}}\circ \sigma)^{-1}$ est bien d\'efinie et est lisse mod\'er\'ee puisque l'inversion et la composition sont des applications lisses mod\'er\'ees dans le contexte fr\'ech\'etique. On en d\'eduit que l'application que nous consid\'erons est bien lisse mod\'er\'ee.$\hfill\square$\\
$\text{}$\\\
Ainsi $\{(\varphi_{\Sigma}(\mathcal{U}_{\Sigma}),\varphi_{\Sigma}^{-1})\,\vert\,\Sigma\in Gr_{k}(M)\}$ est un atlas diff\'erentiable de $Gr_{k}(M)$ et induit canoniquement une topologie $\mathcal{T}$ sur $Gr_{k}(M)\,.$ Montrons que cette topologie est de Hausdorff.
\begin{lemme}
La topologie $\mathcal{T}$ de $Gr_{k}(M)$ est de Hausdorff.
\end{lemme}
\textbf{D\'emonstration.}  Prenons $\Sigma_{1},\Sigma_{2}\in Gr_{k}(M)$ tels que $\Sigma_{1}\neq\Sigma_{2}\,.$ Si ces deux sous-vari\'et\'es orient\'ees se confondent en tant que sous-vari\'et\'es, mais poss\`edent une orientation diff\'erente, alors $\varphi_{\Sigma_{1}}(\mathcal{U}_{\Sigma_{1}})\cap\varphi_{\Sigma_{2}}(\mathcal{U}_{\Sigma_{2}})=\emptyset\,.$ En effet, supposons que $\Sigma$ appartienne \`a cette intersection. Alors $\Sigma$ serait munie de l'orientation induite par le diff\'eomorphisme $\Sigma_{1}\rightarrow\Sigma,\,x\mapsto\tau_{\Sigma_{1}}(s(x))$ o\`u $s$ est une certaine section du fibr\'e normal de $\Sigma_{1}\,.$ De plus, $\Sigma$ serait aussi munie de l'orientation induite par le diff\'eomorphisme $\Sigma_{2}\rightarrow\Sigma,\,x\mapsto\tau_{\Sigma_{2}}(s(x))$, avec $\Sigma_{1}$ et $\Sigma_{2}$ \'etant les m\^emes sous-vari\'et\'es mais orient\'ees diff\'eremment  et $\tau_{\Sigma_{1}}=\tau_{\Sigma_{2}}$, d'o\`u la contradiction.\\\\
Si $\Sigma_{1}\neq\Sigma_{2}$ en tant que sous-vari\'et\'e, il existe $x_{1}\in \Sigma_{1}{\setminus}\Sigma_{2}$ et $\varepsilon>0$ tels que
$$
\exp_{x_{1}}(\text{B}(x_{1},\varepsilon))\cap \tau_{\Sigma_{2}}(\Theta_{\Sigma_{2}})=\emptyset\,,
$$
o\`u $\Theta_{\Sigma_{2}}\subseteq \{v\in N\Sigma_{2}\,\vert\,\Vert x\Vert<\varepsilon\}\,.$  D\`es lors, si l'on choisit $\Theta_{\Sigma_{1}}$ tel que $\Theta_{\Sigma_{1}}\subseteq\{v\in N\Sigma_{1}\,\vert\,\Vert x\Vert<\varepsilon\}$, alors on peut constater que $\varphi_{\Sigma_{1}}(\mathcal{U}_{\Sigma_{1}})\cap\varphi_{\Sigma_{2}}(\mathcal{U}_{\Sigma_{2}})=\emptyset\,.\hfill\square$\\
$\text{}$\\
En r\'esum\'e,
\begin{proposition}[\textit{Hamilton},\,\cite{Hamilton}]
$\text{ }$\:\:L'ensemble $Gr_{k}(M)$ est une vari\'et\'e fr\'ech\'etique lisse mod\'er\'ee et pour $\Sigma\in Gr_{k}(M)$, on a un isomorphisme canonique :
$$
\text{T}_{\Sigma}Gr_{k}(M)\cong\Gamma_{C^{\infty}}(\Sigma, N\Sigma)\,.
$$

\end{proposition}
\begin{remarque}  Tout comme la structure de vari\'et\'e de $C^{\infty}(N,M)$ ne d\'epend pas de la m\'etrique que l'on utilise sur $M$, la structure de vari\'et\'e de $Gr_{k}(M)$ ne d\'epend pas non plus de la m\'etrique $g\,.$
\end{remarque}
\begin{remarque} Notons ${Gr}_{k}^{\vee}(M)$ l'ensemble des sous-vari\'et\'es connexes, compactes, orientables et de dimension $k$ de $M$. Alors, exactement de la m\^eme mani\`ere que pour $Gr_{k}(M)$, on montre que cet ensemble est muni d'une structure de vari\'et\'e mod\'er\'ee et il est clair que $Gr_{k}(M)$ est un rev\^etement \`a deux feuillets de ${Gr}_{k}^{\vee}(M)\,.$
\end{remarque}
\section{L'espace des plongements dans $M$ comme fibr\'e principal sur la Grassmannienne non-lin\'eaire}
Prenons $\Sigma\in Gr_{k}(M)$ et notons $\text{Emb}(\Sigma,M)$ l'espace des plongements de $\Sigma$ dans $M\,.$ Notons aussi
$$
p\,:\,\text{Emb}(\Sigma,M)\rightarrow Gr_{k}(M)\,,
$$
l'application qui est d\'efinie pour $f\in \text{Emb}(\Sigma,M)$ par $p(f):=f(\Sigma)$, cette derni\`ere sous-vari\'et\'e de $M$ \'etant munie de l'orientation naturellement induite par le diff\'eomorphisme $f\,:\,\Sigma\rightarrow f(\Sigma)\,.$
Gr\^ace \`a la proposition suivante et \`a ses corollaires, nous allons montrer que $\text{Emb}(\Sigma,M)$ est une vari\'et\'e fr\'ech\'etique lisse mod\'er\'ee et que l'application $p\,:\,\text{Emb}(\Sigma,M)\rightarrow Gr_{k}(M)$ est lisse. Nous utiliserons pour cela le calcul convenable de Kriegl et Michor (voir \cite{Kriegl-Michor}), car entre des vari\'et\'es fr\'ech\'etiques, une application est lisse au sens de Kriegl-Michor si et seulement si elle est lisse au sens de Hamilton (voir notre succinct appendice)\,. Cela nous am\`enera \`a montrer dans un deuxi\`eme temps que $\text{Emb}(\Sigma,M)$ est l'espace total d'un fibr\'e principal ayant pour base la r\'eunion de certaines composantes connexes de $Gr_{k}(M)\,.$
\begin{proposition}\label{la clef} Soit $E\overset{\pi}{\rightarrow}M$ un fibr\'e vectoriel de rang fini au-dessus de $M$ et $f\in C^{\infty}\big((-\varepsilon,\varepsilon)\times M, E\big)$ une application telle que $f_{0}\,:\,M\rightarrow E,\,x\mapsto f(0,x)$ soit la section nulle de $E\,.$ Alors il existe $\eta>0$ et $\varphi\,:\,(-\eta,\eta)\rightarrow \text{Diff}(M)$ un chemin lisse de $\text{Diff}(M)$ tel que :
\begin{center}
\begin{description}
\item[$(i)$] $f_{t}\circ \varphi_{t}\in \Gamma_{C^{\infty}}(M,E)$ pour tout $t\in (-\eta,\eta)$ \:\:\:(ici $f_{t}(x):=f(t,x)$)\,;\\
\item[$(ii)$] $\varphi_{0}=Id\,.$
\end{description}
\end{center}
\end{proposition}
\textbf{D\'emonstration.} Posons $\psi:\,(-\varepsilon,\,\varepsilon)\times M\rightarrow M,\,(t,x)\mapsto \pi\big( f(t,x)\big)$ et $\psi^{\vee}\,:\,(-\varepsilon,\,\varepsilon)\rightarrow     C^{\infty}(M,M),\,t\mapsto\{M\ni x\mapsto \pi\big( f(t,x)\big)\in M\}\,.$
Etant donn\'e que $\psi=\pi\circ f\,,$ l'application $\psi$ est lisse, ce qui veut exactement dire que $\psi^{\vee}$ est une courbe lisse de $C^{\infty}(M,M)\, $(Voir Appendice, Proposition \ref{caracterisation courbe lisse}). Or, le groupe de Lie $\text{Diff}\,(M)$ \'etant ouvert dans $C^{\infty}(M,M)\,,$ (voir \cite{Kriegl-Michor}, Theorem 43.1.), et puisque $\psi^{\vee}(0)=Id\,,$ on en d\'eduit qu'il existe $\eta>0$ tel que $\psi^{\vee}$ restreint \`a $(-\eta,\,\eta)$ soit une courbe lisse de $\text{Diff}\,(M)\,.$\\\\
Consid\'erons alors le chemin lisse $\varphi\,:\,(-\eta,\,\eta)\rightarrow \text{Diff}\,(M),\,t\mapsto\big(\psi^{\vee}(t)\big)^{-1}\,.$ Pour $x=\big(\psi^{\vee}(t)\big)(y)\in M$, on constate que :
$$
\pi\Big((f_{t}\circ\varphi_{t})(x)\Big)=\pi\Big(f_{t}\Big((\psi_{t}^{\vee})^{-1}(\psi_{t}^{\vee}(y))\Big)\Big)=\pi\Big(f_{t}(y)\Big)=\psi^{\vee}_{t}(y)=x\,.
$$
Donc $\pi\circ(f_{t}\circ\varphi_{t})=Id$ ce qui signifie que $f_{t}\circ\varphi_{t}$ est une section de $E\,.\\\text{}\hfill\square$\\\\
Remarquons que comme corollaire de cette proposition, on retrouve le r\'esultat classique suivant (voir par exemple \cite{Hirsch}, Theorem 1.4, page 37)\,:
\begin{corollaire}
L'ensemble $\text{Emb}(\Sigma,M)$ est ouvert dans $C^{\infty}(\Sigma,M)\,.$ En
particulier, $\text{Emb}(\Sigma,M)$ est naturellement une vari\'et\'e fr\'ech\'etique lisse mod\'er\'ee.
\end{corollaire}
\textbf{D\'emonstration.} Supposons que $\text{Emb}(\Sigma,M)$ ne soit pas ouvert dans $C^{\infty}(\Sigma,M)\,.$ On peut donc trouver une suite $(f_{n})_{n\in\mathbb{N}}$ de $C^{\infty}(\Sigma,M){\setminus}\text{Emb}(\Sigma,M)$ telle que $f_{n}\rightarrow f$ pour un certain $f\in \text{Emb}(\Sigma,M)\,.$ Soit alors $(\mathcal{U}_{f},\varphi_{f})$ une carte de $C^{\infty}(\Sigma,M)$ centr\'ee en $f$ et telle que $\varphi_{f}(\mathcal{U}_{f})$ soit convexe. Rappelons que l'on peut construire la carte $(\mathcal{U}_{f},\varphi_{f})$ en prenons $\varphi_{f}(\mathcal{U}_{f})$ un voisinage de 0 suffisamment petit de l'espace des sections $\Gamma_{C^{\infty}}(\Sigma,\,f^{*}TM)$ et $\varphi_{f}^{-1}:\,\varphi_{f}(\mathcal{U}_{f})\rightarrow \mathcal{U}_{f}\subseteq C^{\infty}(\Sigma,M)\,,$ l'application qui est d\'efinie par $\varphi^{-1}_{f}(X)(x):=\text{exp}_{f(x)}\,(X_{x})$ pour $X\in \Gamma_{C^{\infty}}(\Sigma,\,f^{*}TM)$ et $x\in \Sigma\,.$ Puisque $f_{n}\rightarrow f$, nous pouvons supposer que $f_{n}\in \mathcal{U}_{f}$ pour tout $n\in\mathbb{N}$, et ainsi consid\'erer la suite de sections $\big(\varphi_{f}(f_{n})\big)_{n\in \mathbb{N}}$ de $\varphi_{f}(\mathcal{U}_{f})\subseteq\Gamma_{C^{\infty}}(\Sigma,f^{*}TM)\,.$ Or, $\Gamma_{C^{\infty}}(\Sigma,f^{*}TM)$ \'etant un espace de Fr\'echet, nous pouvons supposer (voir Appendice, Proposition \ref{special curve}) qu'il existe une courbe lisse de sections $s\,:\,\mathbb{R}\rightarrow\Gamma_{C^{\infty}}(\Sigma,f^{*}TM)$ telle que :
$$
s_{0}=s(0)=\varphi_{f}(f)\,\,\,\,et\,\,\,\,s\bigg(\dfrac{1}{n}\bigg)=\varphi_{f}(f_{n})
$$
pour tout $n\in \mathbb{N}\,.$ Si l'on construit $s$ de la m\^eme mani\`ere que dans le ``special curve lemma'' de \cite{Kriegl-Michor}, on constate que $s$ est \`a valeurs dans $\varphi_{f}(\mathcal{U}_{f})\,.$ En effet, $s(\Sigma)$ est le polygone d'arr\^etes les $\varphi_{f}(f_{n})$ et l'on a choisit $\varphi_{f}(\mathcal{U}_{f})$ convexe. Notons alors
$$
g\,:\,\mathbb{R}\times\Sigma\rightarrow M,\,\,(t,x)\mapsto\varphi_{f}^{-1}(s_{t})(x)\,.
$$
Par construction on a pour tout $n\in \mathbb{N}\,:$
$$
g_{0}=f\,\,\,\,\text{et}\,\,\,\,g_{\frac{1}{n}}=f_{n}\,.
$$
Notons $W:=f(\Sigma)=g_{0}(\Sigma)\,.$ Pour $t$ assez petit, $g_{t}(\Sigma)\subseteq\tau_{W}(\Theta_{W})$ et nous pouvons d\`es lors consid\'erer l'application $W\ni x\rightarrow(\tau_{W}^{-1}\circ g_{t}\circ g_{0}^{-1})(x)\in NW\,.$ Cette derni\`ere application v\'erifie les hypoth\`eses de la Proposition \ref{la clef}, il existe donc une courbe $\varphi_{t}$ de $\text{Diff}(W)$ telle que :
$$
\sigma_{t}\,:\,x\in W\mapsto(\tau_{W}^{-1}\circ g_{t}\circ g_{0}^{-1}\circ\varphi_{t})(x)\in NW
$$
soit une section du fibr\'e normal de $W\,.$ Mais alors, pour $n$ assez grand,
$$
f_{n}=g_{\frac{1}{n}}=\tau_{W}\circ\sigma_{\frac{1}{n}}\circ\varphi_{\frac{1}{n}}^{-1}\circ g_{0}
$$
est un plongement de $\Sigma$ dans $M$ ce qui est une contradiction. Ainsi, $\text{Emb}(\Sigma,M)$ est bien un ouvert de $C^{\infty}(\Sigma,M)\,.\hfill\square$
\begin{corollaire}
L'application $p\,:\,\text{Emb}(\Sigma,M)\rightarrow Gr_{k}(M)$ est lisse et pour un chemin lisse $f_{t}$ de $\text{Emb}(\Sigma,M)$, on a la formule :
$$
\dfrac{d}{dt}\bigg\vert_{t_{0}}\,p(f_{t})=\bigg(\dfrac{\partial f}{\partial t}(t_{0})\bigg)^{\bot}
$$
o\`u $\big(\frac{\partial f}{\partial t}(t_{0})\big)^{\bot}$ est la section de $\Gamma_{C^{\infty}}(f_{t_{0}}(\Sigma),Nf_{t_{0}}(\Sigma))$ qui est d\'efinie pour $x\in \Sigma$ par
$\Big(\dfrac{\partial f}{\partial t}(t_{0})\Big)^{\bot}_{f_{t_{0}}(x)}:=pr\big(\frac{\partial f}{\partial t}(t_{0},x)\big)
\,,$ $pr$ \'etant la projection orthogonale sur le fibr\'e normal de $f_{t_{0}}(\Sigma)\,.$
\end{corollaire}
\textbf{D\'emonstration.} Prenons $f_{t}$ un chemin lisse de $\text{Emb}(\Sigma,M)\,.$ Nous devons montrer que $p(f_{t})$ est un chemin lisse de $Gr_{k}(M)$ afin de v\'erifier la lissit\'e de $p$ au sens de Kriegl-Michor. Fixons $t_{0}\in\mathbb{R}$ et notons $W:=p(f_{t_{0}})\,.$ Pour $\varepsilon>0$  suffisamment petit, l'application $(t,x)\in (t_{0}-\varepsilon,t_{0}+\varepsilon)\times W \mapsto (\tau_{W}^{-1}\circ f_{t}\circ f_{t_{0}}^{-1})(x)\in NW$, satisfait les hypoth\`eses de la Proposition \ref{la clef}\,. Il existe donc un chemin lisse $\varphi_{t}$ de diff\'eomorphismes de $W$ tel que
$$
x\in W\mapsto(\tau_{W}^{-1}\circ f_{t}\circ f_{t_{0}}^{-1}\circ\varphi_{t})(x)\in NW
$$
soit une section de $\Gamma_{C^{\infty}}(W,NW)$ d\`es que $t$ est suffisamment petit. Mais ceci nous donne justement la possibilit\'e d'exprimer $p(f_{t})$ au voisinage de $t_{0}$ dans la carte $\big(\varphi_{W}(\mathcal{U}_{W}),\varphi_{W}^{-1}\big)$ :
$$
\varphi_{W}^{-1}\big(p(f_{t})\big)=\tau_{W}^{-1}\circ f_{t}\circ f_{t_{0}}^{-1}\circ\varphi_{t}\in \Gamma_{C^{\infty}}(W,NW)\,.
$$
Cette derni\`ere courbe de sections \'etant lisse, il en r\'esulte que $p$ est lisse.\\\\
Pour la formule de la diff\'erentielle de $p$, notons $W:=f_{t_{0}}(\Sigma)\,.$ En identifiant $\text{T}_{W}Gr_{k}(M)$ \`a $\Gamma_{C^{\infty}}(W,NW)$, on a :
$$
\dfrac{d}{dt}\bigg\vert_{t_{0}}\,p(f_{t})=\dfrac{d}{dt}\bigg\vert_{t_{0}}\,\varphi_{W}^{-1}\big(p(f_{t})\big)=\dfrac{d}{dt}\bigg\vert_{t_{0}}\,\big(\tau_{W}^{-1}\circ f_{t}\circ f_{t_{0}}^{-1}\circ\varphi_{t}\big)
$$
et pour $x\in W,$
\begin{eqnarray*}
\dfrac{d}{dt}\bigg\vert_{t_{0}}\,\big(\tau_{W}^{-1}\circ f_{t}\circ f_{t_{0}}^{-1}\circ\varphi_{t}\big)(x)&=&(\tau_{W}^{-1})_{*_{x}}
\dfrac{d}{dt}\bigg\vert_{t_{0}}\,\big( f_{t}\circ f_{t_{0}}^{-1}\circ\varphi_{t}\big)(x)\\
&=&(\tau_{W}^{-1})_{*_{x}}\bigg[\underbrace{\dfrac{\partial f}{\partial t}(t_{0},f_{t_{0}}^{-1}(x))}_{\in \text{T}_{x}M}+\underbrace{\dfrac{\partial \varphi}{\partial t}(t_{0},x)}_{\in \text{T}_{x}W}\bigg]\,.
\end{eqnarray*}
Par construction de $\varphi_{t}$, il vient :
$$
\dfrac{\partial f}{\partial t}(t_{0},f_{t_{0}}^{-1}(x))+\dfrac{\partial \varphi}{\partial t}(t_{0},x)\in N_{x}W
$$
ce qui implique, puisque $\dfrac{\partial \varphi}{\partial t}(t_{0},x)\in \text{T}_{x}W$, que
$$
\dfrac{\partial f}{\partial t}(t_{0},f_{t_{0}}^{-1}(x))
+\dfrac{\partial \varphi}{\partial t}(t_{0},x)=\bigg(\dfrac{\partial f}{\partial t}(t_{0})\bigg)_{x}^{\perp}\,.
$$
Ainsi,
\begin{eqnarray*}
&&\dfrac{d}{dt}\bigg\vert_{t_{0}}\,\big(\tau_{W}^{-1}\circ f_{t}\circ f_{t_{0}}^{-1}\circ\varphi_{t}\big)(x)=(\tau_{W}^{-1})_{*_{x}}\bigg(\dfrac{\partial f}{\partial t}(t_{0})\bigg)_{x}^{\perp}\\
&=&\dfrac{d}{du}\bigg\vert_{0}\,\tau_{W}^{-1}\bigg[exp_{x}\bigg(u\bigg(\dfrac{\partial f}{\partial t}(t_{0})\bigg)_{x}^{\perp}\bigg)\bigg]
=\dfrac{d}{du}\bigg\vert_{0}\,(\tau_{W}^{-1}\circ\tau_{W})\bigg(x,u\bigg(\dfrac{\partial f}{\partial t}(t_{0})\bigg)_{x}^{\perp}\bigg)\\
&=&\dfrac{d}{du}\bigg\vert_{0}\,\bigg(x,u\bigg(\dfrac{\partial f}{\partial t}(t_{0})\bigg)_{x}^{\perp}\bigg)
=\bigg(\dfrac{\partial f}{\partial t}(t_{0})\bigg)_{x}^{\perp}\,.
\end{eqnarray*}
Par suite,
$$
\dfrac{d}{dt}\bigg\vert_{t_{0}}\,p(f_{t})=\bigg(\dfrac{\partial f}{\partial t}(t_{0})\bigg)^{\perp}\,.
$$
$\text{}\hfill\square$\\\\
\begin{remarque} Au vu de la formule de la diff\'erentielle de $p$, il semblerait que cette derni\`ere application d\'epende ``plus'' que de la structure diff\'erentiable de $Gr_{k}(M)$ puisque la m\'etrique utilis\'ee apparait dans la formule de la diff\'erentielle de $p\,.$ En fait, il ne faut pas oublier que pour $W\in Gr_{k}(M)$, $\text{T}_{W}Gr_{k}(M)$ n'est pas \'egal \`a l'espace $\Gamma_{C^{\infty}}(W,NW)$ mais lui est seulement isomorphe via une r\'ealisation n\'ecessitant la m\'etrique $g\,.$
\end{remarque}
A pr\'esent, afin de pouvoir consid\'erer certains fibr\'es principaux, introduisons, pour $\Sigma\in Gr_{k}(M)$, les notations suivantes :
\begin{description}
\item[$(i)$] $p\,:\text{Emb}(\Sigma,M)\,\rightarrow Gr_{k}(M),\,\,f\mapsto p(f)$ la projection canonique;\\
\item[$(ii)$] $p_{\Sigma}\,:\text{Emb}(\Sigma,M)\,\rightarrow Gr(\Sigma,M):=p(\text{Emb}(\Sigma,M)),\,\,f\mapsto p(f)$\,;\\
\item[$(iii)$] $\lambda\,:\,\text{Emb}(\Sigma,M)\times \text{Diff}^{+}(\Sigma)\rightarrow \text{Emb}(\Sigma,M),\,\,(f,\varphi)\mapsto f\circ\varphi$ l'action naturelle \`a droite de $\text{Diff}^{+}(\Sigma)$ sur $\text{Emb}(\Sigma,M)\,.$\\
\end{description}
Nous allons montrer par une s\'erie de lemmes que $\text{Emb}(\Sigma,M)$ est un $\text{Diff}^{+}(\Sigma)$-fibr\'e principal au-dessus de $Gr(\Sigma,M)\,.$
\begin{lemme}\label{composition}
Soient $U,V$ deux ouverts de $M$ d'intersection non nulle et  $\Sigma_{0},\Sigma_{1},W$
trois sous-vari\'et\'es de $M$ telles que :
$$
\Sigma_{0}\subseteq U,\,\,\,\Sigma_{1}\subseteq V\,\,\,\,\,\,\,\text{et}\,\,\,\,\,\,\,W\subseteq U\cap V\,.
$$
Soient aussi $\beta$ un chemin continu de $\text{Emb}(\Sigma_{0}, U)$ et $\tilde{\beta}$ un chemin continu de $\text{Emb}(W, V)$ tels que :
$$
\beta(0)=j_{\Sigma_{0}},\,\,\beta(1)(\Sigma_{0})=W,\,\,\tilde{\beta}(0)=j_{W}\,\,\,\,\text{et}\,\,\,\,\tilde{\beta}(1)(W)=\Sigma_{1}
$$
o\`u $j_{\Sigma_{0}}\,:\,\Sigma_{0}\hookrightarrow M$ et $j_{W}\,:\,W\hookrightarrow M$ sont les inclusions canoniques.
Alors, l'application $\gamma\,:\,[0,2]\rightarrow \text{Emb}(\Sigma_{0},U\cup V)$ d\'efinie par :
$$
\gamma(t)=
\left\lbrace
\begin{array}{c}
 \beta(t)\,\,\,\,\text{pour}\,\,\,\,t\in [0,1]\,;\\
\tilde{\beta}(t-1)\circ\beta(1)\,\,\,\,\text{pour}\,\,\,\,t\in[1,2],
\end{array}
\right.
$$
est un chemin continu de $\text{Emb}(\Sigma_{0}, U\cup V)\,.$
\end{lemme}
\textbf{D\'emonstration.} Consid\'erons l'application
$$
\vartheta\,:\,\text{Emb}(W, V)\rightarrow\text{Emb}(\Sigma_{0},U\cup V),\,\,\,\rho\mapsto\rho\circ\beta(1)\,.
$$
En utilisant les courbes lisses de $\text{Emb}(W, V)$, il est imm\'ediat que $\vartheta$ est une application lisse, en particulier, $\vartheta$ est continue. Mais alors :\begin{description}
\item[$(i)$]  $\gamma$ est clairement continue sur $[0,1]\,;$\\
\item[$(ii)$]  $\gamma(t)=(\vartheta\circ\tilde{\beta})(t-1)$ pour $t\in [1,2]$ et donc $\gamma$ est continue sur [1,2].
\end{description}
Il en r\'esulte que $\gamma\,:\,[0,2]\rightarrow \text{Emb}(\Sigma_{0},U\cup V)$ est bien un chemin continu de $ \text{Emb}(\Sigma_{0},U\cup V)\,.\hfill\square$
\begin{lemme}\label{lift} Soient $\Sigma_{0},\Sigma_{1}$ deux \'el\'ements de $Gr_{k}(M)$ et $\alpha\,:\,[0,1]\rightarrow Gr_{k}(M)$ un chemin continu tel que $\alpha(0)=\Sigma_{0}$ et $\alpha(1)=\Sigma_{1}\,.$ Alors il existe $\beta\,:\,[0,1]\rightarrow\text{Emb}(\Sigma_{0},M)$, un chemin continu de $\text{Emb}(\Sigma_{0},M)$ tel que :
$$
\beta(0)=j_{\Sigma_{0}}\,,\,\,\,\,(p\circ\beta)(0)=\alpha(0)=\Sigma_{0}\,\,\,\,\text{et}\,\,\,\,(p\circ\beta)(1)=\alpha(1)=\Sigma_{1}
$$
o\`u $p\,:\, \text{Emb}(\Sigma_{0},M)\rightarrow Gr_{k}(M)$ est la projection canonique et $j_{\Sigma_{0}}\,:\,\Sigma_{0}\hookrightarrow M$ l'inclusion canonique.
\end{lemme}
\textbf{D\'emonstration.} Puisque $\alpha$ est continu, l'ensemble $\alpha([0,1])$ est compact et peut donc \^etre recouvert par un nombre fini de carte. Pour simplifier, supposons que
$$
\alpha([0,1])\subseteq\varphi_{\Sigma_{0}}(\mathcal{U}_{\Sigma_{0}})\cup\varphi_{\Sigma_{1}}(\mathcal{U}_{\Sigma_{1}})\,,
$$
 les notations \'etant celles pr\'ec\'edement introduites. Prenons $W$ un \'el\'ement de $\varphi_{\Sigma_{0}}(\mathcal{U}_{\Sigma_{0}})\cap\varphi_{\Sigma_{1}}(\mathcal{U}_{\Sigma_{1}})\,.$ On a :
$$
W=\varphi_{\Sigma_{0}}(s_{0})\,\,\,\,\text{et}\,\,\,\,W=\varphi_{\Sigma_{1}}(s_{1})
$$
pour un certain $s_{0}\in \mathcal{U}_{\Sigma_{0}}$ et un certain $s_{1}\in \mathcal{U}_{\Sigma_{1}}\,.$
On peut alors consid\'erer les applications :
$$
\beta\,:\,[0,1]\rightarrow\text{Emb}(\Sigma_{0},\tau_{\Sigma_{0}}(\Theta_{\Sigma_{0}})),\,\,t\mapsto\{x\in\Sigma_{0}\mapsto\text{exp}_{x}(ts_{0}(x))\}
$$
et
$$
\tilde{\beta}\,:\,[0,1]\rightarrow\text{Emb}(W,\tau_{\Sigma_{1}}(\Theta_{\Sigma_{1}})),\,\,t\mapsto\{x\in W\mapsto\text{exp}_{\tilde{p}(x)}((1-t)s_{1}(\tilde{p}(x)))\}\,.
$$
o\`u $\tilde{p}\,:\,\tau_{\Sigma_{1}}(\Theta_{\Sigma_{1}})\rightarrow\Sigma_{1},\,\,\tau_{\Sigma_{1}}\big((x,v)\big)\mapsto x$ pour $(x,v)\in N\Sigma_{1}$ est la projection canonique. Ces deux applications, $\beta$ et $\tilde{\beta}\,,$ sont manifestement continues puisque l'on peut les \'etendre en des courbes lisses. Il en r\'esulte par le Lemme \ref{composition} (et apr\`es reparam\'etrage), qu'il existe une courbe continue $\gamma\,:\,[0,1]\rightarrow\text{Emb}(\Sigma_{0},\tau_{\Sigma_{0}}(\Theta_{\Sigma_{0}})\cup\tau_{\Sigma_{1}}(\Theta_{\Sigma_{1}}))\subseteq\text{Emb}(\Sigma_{0},M)$ telle que $\gamma(0)=\beta(0)=j_{\Sigma_{0}}$ et $\gamma(1)=\tilde{\beta}(1)\circ \beta(1)\,.$ D'o\`u :
$$
(p\circ\gamma)(0)=p\big(\gamma(0)\big)=p(j_{\Sigma_{0}})=\Sigma_{0}\,.
$$
De plus, puisque
$$
\gamma(1)(\Sigma_{0})=\big(\tilde{\beta}(1)\circ \beta(1)\big)(\Sigma_{0})=\tilde{\beta}(1)(W)=\Sigma_{1}\,,
$$
il suffit pour montrer que $p(\gamma(1))=\Sigma_{1}$, de v\'erifier que $[\gamma(1)^{*}\mu_{1}]=[\mu_{0}]$ o\`u $[\mu_{i}]$ est l'orientation de $\Sigma_{i}$ (i=0,1). Mais cela d\'ecoule de la d\'efinition m\^eme des cartes $\varphi_{\Sigma_{i}}(\mathcal{U}_{\Sigma_{i}})\,.$  En effet, d'apr\`es cette d\'efinition, l'orientation de $W$ est donn\'ee \`a la fois par $[(\beta(1)^{-1})^{*}\mu_{0}]$ et par $[\tilde{\beta}(1)^{*}\mu_{1}]\,,$ et comme $\gamma(1)=\tilde{\beta}(1)\circ \beta(1)\,,$ on en d\'eduit que $p(\gamma(1))=(\Sigma_{1},[\mu_{1}])=\Sigma_{1}\,.\hfill\square$\\\\
\begin{corollaire}
L'ensemble $Gr(\Sigma,M)$ est une r\'eunion de composantes connexes de $Gr_{k}(M)\,.$ En particulier, $Gr(\Sigma,M)$ est une vari\'et\'e mod\'er\'ee.
\end{corollaire}
\textbf{D\'emonstration.}  Soient $f\in \text{Emb}(\Sigma,M)$ et $\Sigma_{1}\in Gr_{k}(M)$ un \'el\'ement appartenant \`a la m\^eme composante connexe dans $Gr_{k}(M)$ que $p_{\Sigma}(f)=:\Sigma_{0}\,.$ Pour montrer le lemme, il suffit de montrer que $\Sigma_{1}\in Gr(\Sigma,M)\,.$\\\\
Prenons $\alpha\,:\,[0,1]\rightarrow Gr_{k}(M)$, un chemin continu de $Gr_{k}(M)$ tel que :
$$
\alpha(0)=p_{\Sigma}(f)=\Sigma_{0}\:\:\:\:\text{et}\:\:\:\:\alpha(1)=\Sigma_{1}\,.
$$
D'apr\`es le Lemme \ref{lift}, il existe un chemin $\beta\,:\,[0,1]\rightarrow\text{Emb}(\Sigma_{0},M)$ tel que :
$$
\beta(0)=j_{\Sigma_{0}}\,,\,\,\,\,(p_{\Sigma_{0}}\circ\beta)(0)=\alpha(0)=\Sigma_{0}\,\,\,\,\text{et}\,\,\,\,(p_{\Sigma_{0}}\circ\beta)(1)=\alpha(1)=\Sigma_{1}\,.
$$
Si l'on regarde $f$ comme une application \`a valeurs dans $\Sigma_{0}=p_{\Sigma}(f)$, alors on constate que
$$
p_{\Sigma}(\beta(1)\circ f)=\Sigma_{1}\:\:\:\:\text{avec}\:\:\:\:\:\beta(1)\circ f\in \text{Emb}(\Sigma,M)\,.
$$
Ainsi, $\Sigma_{1}\in Gr(\Sigma,M)\,.\hfill\square$
\begin{lemme}\label{section locale}
L'application $p_{\Sigma}\,:\,\text{Emb}(\Sigma,M)\rightarrow Gr(\Sigma,M)$ admet des sections locales.
\end{lemme}
\textbf{D\'emonstration.} Soit $W\in Gr(\Sigma,M)$ et soit $f\in \text{Emb}(\Sigma,M)$ tel que $p_{\Sigma}(f)=W\,.$ On peut constater que l'application $\sigma \,:\,\varphi_{W}(\mathcal{U}_{W})\rightarrow \text{Emb}(\Sigma,M)$ d\'efinie par
$$
\sigma(\varphi_{W}(s))(x):=exp_{f(x)}\,s(f(x))
$$
est une section locale de $p_{\Sigma}\,.\hfill\square$
\begin{lemme}
L'action $\lambda\,:\,\text{Emb}(\Sigma,M)\times \text{Diff}^{+}(M)\rightarrow \text{Emb}(\Sigma,M)$ est libre. De plus, pour $W\in Gr(\Sigma,M)$ et $f\in \text{Emb}(\Sigma,M)$ telle que $p_{\Sigma}(f)=W\,,$ on a $p_{\Sigma}^{-1}(W)=\mathcal{O}_{f}$ o\`u $\mathcal{O}_{f}$ est l'orbite de $f$ pour l'action $\lambda\,.$
\end{lemme}
\textbf{D\'emonstration.} La libert\'e de $\lambda$ est \'evidente. Montrons que $p_{\Sigma}^{-1}(W)=\mathcal{O}_{f}\,.$ Notons $[\mu]$ l'orientation de $\Sigma\,.$ Par la suite, nous noterons $f^{-1}$ l'unique application lisse de $W$ dans $\Sigma$ v\'erifiant $f^{-1}\circ f=id_{\Sigma}$ (et de m\^eme pour $g$)\,. On a :
\begin{eqnarray*}
g\in p_{\Sigma}^{-1}(W)&\Leftrightarrow& g(\Sigma)=f(\Sigma)\:\:\:\text{et}\:\:\:[(g^{-1})^{*}\mu]=[(f^{-1})^{*}\mu]\\
&\Leftrightarrow& g=f\circ \varphi\:\:\text{avec}\:\:\varphi=f^{-1}\circ g \in \text{Diff}(\Sigma)\\
&& \text{et}\:\:\:[(g^{-1})^{*}\mu]=[(f^{-1})^{*}\mu]\\
&\Leftrightarrow&g=f\circ\varphi\:\:\text{avec}\:\:\varphi\in\text{Diff}^{+}(\Sigma)\\
&\Leftrightarrow& g=\lambda(f,\varphi)\:\:\text{avec}\:\:\varphi\in\text{Diff}^{+}(\Sigma)\,.
\end{eqnarray*}
Ainsi, $p_{\Sigma}^{-1}(p_{\Sigma}(f))=\mathcal{O}_{f}\,.\text{}\hfill\square$\\\\
\begin{lemme}\label{bien lisse}
Pour $W\in Gr(\Sigma,M)$ et $f\in \text{Emb}(\Sigma,M)$ telle que $p_{\Sigma}(f)=\Sigma$\,, l'application
$$
\Lambda\::\:p_{\Sigma}^{-1}(\varphi_{W}(\mathcal{U}_{W}))\longrightarrow\text{Diff}^{+}(\Sigma),\,\,\,g\mapsto\sigma\Big(p_{\Sigma}(g)\Big)^{-1}\circ g
$$
est lisse mod\'er\'ee (ici $\sigma$ correspond \`a la section construite dans le Lemme \ref{section locale}).
\end{lemme}
\textbf{D\'emonstration.}  Remarquons que $\Lambda$ est bien d\'efinie et est l'unique application v\'erifiant $\lambda\Big(\sigma\big(p_{\Sigma}(g)\big),\,\Lambda(g)\Big)=g$ pour tout $g\in p_{\Sigma}^{-1}(\varphi_{W}(\mathcal{U}_{W}))\,.$\\\\
Montrons que $\Lambda$ est lisse mod\'er\'ee. Pour $g\in p_{\Sigma}^{-1}(\varphi_{W}(\mathcal{U}_{W}))$ et $x\in \Sigma$ on a :
\begin{eqnarray*}
\Lambda(g)(x)=\Big(\Big(\sigma\big(p_{\Sigma}(g)\big)\Big)^{-1}\circ g\Big)(x)&\Rightarrow&\sigma\big(p_{\Sigma}(g)\big)\big(\Lambda(g)(x)\big)=g(x)\,.
\end{eqnarray*}
Notons $s\in \Gamma_{C^{\infty}}(W, NW)$ l'unique section de $\Gamma_{C^{\infty}}(W, NW)$ v\'erifiant $p_{\Sigma}(g)=\varphi_{W}(s)\,.$ On a alors :
\begin{eqnarray*}
\sigma\big(\varphi_{W}(s)\big)\big(\Lambda(g)(x)\big)=g(x)&\Rightarrow&\text{exp}_{(f\circ\Lambda(g))(x)}\,\,\Big(\big(s\circ f\circ \Lambda(g)\big)(x)\Big)=g(x)\\
&\Rightarrow&\big(s\circ f\circ \Lambda(g)\big)(x)=\tau_{W}^{-1}\big(g(x)\big)\\
&\Rightarrow&\big( f\circ \Lambda(g)\big)(x)=\big(\pi_{NW}\circ\tau_{W}^{-1}\circ g\big)(x)\\
&\Rightarrow& \Lambda(g)(x)=\big(f^{-1}\circ\pi_{NW}\circ\tau_{W}^{-1}\circ g\big)(x)\,.\\
\end{eqnarray*}
Ainsi, $\Lambda(g)=f^{-1}\circ\pi_{NW}\circ\tau_{W}^{-1}\circ g\,,$ et l'on peut remarquer que l'application $f$ \'etant fix\'ee,  $f^{-1}$ est une application lisse ind\'ependante de $g\in p_{\Sigma}^{-1}(\varphi_{W}(\mathcal{U}_{W}))\,.$ On en d\'eduit que $\Lambda$ est bien une application lisse mod\'er\'ee.\\$\text{}\hfill\square$\\\\
De cette succession de lemmes, on en d\'eduit :
\begin{theoreme}
L'application $p_{\Sigma}\,:\,\text{Emb}(\Sigma,M)\rightarrow Gr(\Sigma,M)$ est un $\text{Diff}^{+}(\Sigma)$-fibr\'e principal mod\'er\'e pour l'action $\lambda\,.$
\end{theoreme}
\textbf{D\'emonstration.}  Prenons $W\in \text{Gr}(\Sigma,\,M)$ et choisissons $f\in \text{Emb}(\Sigma,M)$ telle que $p_{\Sigma}(f)=W\,.$ On peut consid\'erer le diagramme commutatif suivant :
\begin{center}
$
\begin{diagram}
\node{p_{\Sigma}^{-1}\big(\varphi_{W}(\mathcal{U}_{W})\big)} \arrow[2]{e,t}{\displaystyle\Psi} \arrow{se,b}{\displaystyle p_{\Sigma}} \node[2]{\varphi_{W}(\mathcal{U}_{W})\times\text{Diff}^{+}(\Sigma)} \arrow{sw,b}{\displaystyle pr_{1}}\\
\node[2]{\varphi_{W}(\mathcal{U}_{W})}
\end{diagram}
$
\end{center}
o\`u $\Psi(g):=(p_{\Sigma}(g),\,\Lambda(g))\,.$ \\
D'apr\`es le Lemme \ref{bien lisse}, $\Psi$ est une application lisse mod\'er\'ee. Cette application est de plus $\text{Diff}^{+}$($\Sigma$)-\'equivariante, d'inverse lisse mod\'er\'ee $\Psi^{-1}\big(\varphi_{W}(s),\,\varphi\big)=\sigma\big(\varphi_{W}(s)\big)\circ\varphi\,.$ On construit ainsi des trivialisations de $\text{Emb}(\Sigma,\,M)$ faisant de $\text{Emb}(\Sigma,\,M)$ un $\text{Diff}^{+}$($\Sigma$)-fibr\'e principal au-dessus de $\text{Gr}(\Sigma,\,M)\,.\hfill\square$

\section{Homog\'en\'eit\'e des composantes connexes de $Gr_{k}(M)$ sous l'action des diff\'eomorphismes de $M$}
Pour $\Sigma\in Gr_{k}(M)$, nous savons que la composante connexe $\big(Gr_{k}(M)\big)_{\Sigma}$ de $ Gr_{k}(M)$ contenant $\Sigma$ est connexe et localement connexe par arcs (l'espace mod\`ele \'etant de Fr\'echet) et donc, $\big(Gr_{k}(M)\big)_{\Sigma}$ est aussi connexe par arcs. On a alors, tout comme en dimension finie :
\begin{proposition}\label{connexite}
La composante connexe $\big(Gr_{k}(M)\big)_{\Sigma}$ est connexe par arcs pour des arcs lisses.
\end{proposition}
Pour montrer ce r\'esultat, nous avons besoin d'un lemme que l'on peut d\'eduire de \cite{Hirsch} (voir exercice 3.b, section 8.1, page 182 de \cite{Hirsch}).
\begin{lemme}\label{Hirsch}
Si $\beta\,:\,[0,1]\rightarrow\text{Emb}(\Sigma,M)$ est un chemin continu, alors il existe une application lisse $F\,:\, [0,1]\times\Sigma\rightarrow M$ telle que :
\begin{center}
\begin{description}
\item[$(i)$]  l'application $F_{t}\,:\,\Sigma\rightarrow M,\,\,x\mapsto F(t,x)$ soit un plongement pour tout $t\in [0,1]$\,;\\
\item[$(ii)$]  $F_{0}(\Sigma)=\beta(0)(\Sigma)$ et $F_{1}(\Sigma)=\beta(1)(\Sigma)\,.$
\end{description}
\end{center}
\end{lemme}
\textbf{D\'emonstration de la Proposition \ref{connexite}.}  Prenons $\alpha\,:\,[0,1]\rightarrow\big(Gr_{k}(M)\big)_{\Sigma}$ un chemin continu de $\big(Gr_{k}(M)\big)_{\Sigma}\,.$ Notons $\Sigma_{0}:=\alpha(0)$ et $\Sigma_{1}:=\alpha(1)\,.$ D'apr\`es le Lemme \ref{lift}, il existe $\beta\,:[0,1]\rightarrow\text{Emb}(\Sigma_{0},M)$ un chemin continu de $\text{Emb}(\Sigma_{0},M)$ tel que :
$$
(p\circ\beta)(0)=\alpha(0)\,\,\,\,\text{et}\,\,\,\,(p\circ\beta)(1)=\alpha(1)\,.
$$
Mais alors, d'apr\`es le Lemme \ref{Hirsch}, nous pouvons trouver $\varepsilon>0$ et une application lisse $F\,:\,]-\varepsilon,1+\varepsilon[\times\Sigma_{0}\rightarrow M$ telle que :
\begin{center}
\begin{description}
\item[$(i)$]  l'application $F_{t}\,:\,\Sigma_{0}\rightarrow M,\,x\mapsto F(t,x)$ soit un plongement pour tout $t\in [0,1]\,;$\\
\item[$(ii)$]  $F_{0}(\Sigma_{0})=\beta(0)(\Sigma_{0})$ et $F_{1}(\Sigma_{0})=\beta(1)(\Sigma_{0})\,.$
\end{description}
\end{center}
Il en r\'esulte que l'application $t\in ]-\varepsilon,1+\varepsilon[\rightarrow\text{Emb}(\Sigma_{0},M),\,t\mapsto F_{t}$ est une courbe lisse de $\text{Emb}(\Sigma_{0},M)$ pour $\varepsilon$ suffisamment petit (car $\text{Emb}(\Sigma_{0},M)$ est ouvert dans $C^{\infty}(\Sigma_{0},M)$).\\ Par suite, $p\circ F_{t}$ est une courbe lisse de $\big(Gr_{k}(M)\big)_{\Sigma}$ v\'erifiant :
$$
p\circ F_{0}=F_{0}(\Sigma_{0})=\beta(0)(\Sigma_{0})=\alpha(0)=\Sigma_{0}
$$
$$
\text{et}\,\,\,\,p\circ F_{1}=F_{1}(\Sigma_{0})=\beta(1)(\Sigma_{0})=\alpha(1)=\Sigma_{1}
$$
ce qui montre la proposition.$\hfill\square$\\\\\\

A pr\'esent, consid\'erons $\text{Diff}^{\,0}(M)$, la composante connexe de $\text{Diff}(M)$ contenant l'\'el\'ement neutre $Id_{M}$ ainsi que son action naturelle sur $\big(Gr_{k}(M)\big)_{\Sigma}$ :
$$
\vartheta\,:\,\text{Diff}^{\,0}(M)\times\big(Gr_{k}(M)\big)_{\Sigma}\rightarrow\big(Gr_{k}(M)\big)_{\Sigma}\,,\,\,(\varphi,W)\rightarrow \varphi(W)\,.
$$
On a alors le r\'esultat d'homog\'en\'eit\'e suivant:
\begin{theoreme}\label{theoreme 2}
L'action de $\text{Diff}^{\:\,0}(M)$ sur $\big(Gr_{k}(M)\big)_{\Sigma}$ est transitive.
\end{theoreme}
\textbf{D\'emonstration.}   $\text{}$\:\:\:\:Soient $\Sigma_{0}$ et $\Sigma_{1}$ deux \'el\'ements de $\big(Gr_{k}(M)\big)_{\Sigma}$ et $\alpha\,:\,[0,1]\rightarrow \big(Gr_{k}(M)\big)_{\Sigma}$ une courbe continue joignant $\Sigma_{0}$ et $\Sigma_{1}\,.$ Tout comme dans la d\'emonstration de la Proposition \ref{connexite}, nous pouvons trouver une application lisse $F\,:\,[0,1]\times\Sigma_{0}\rightarrow M$ telle que :
$$
F_{0}(\Sigma_{0})=\Sigma_{0}\,\,\,\,\text{et}\,\,\,\, F_{1}(\Sigma_{0})=\Sigma_{1}\,
$$
et telle que $F_{t}$ soit un plongement pour tout $t\in[0,1]\,.$
Mais alors, d'apr\`es un r\'esultat classique de topologie diff\'erentielle (voir Th\'eor\`eme 1.3, chapitre 8, page 180 de \cite{Hirsch}), nous pouvons trouver une application lisse $\widetilde{F}\,:\,[0,1]\times M\rightarrow M$ v\'erifiant  pour tout $t\in [0,1]$:
\begin{center}
\begin{description}
\item[$(i)$] $\widetilde{F}_{t}\in \text{Diff}(M)\,;$\\
\item[$(ii)$] $\widetilde{F}_{0}=Id$ et $F_{t}=\widetilde{F}_{t}\vert_{\Sigma_{0}}\,.$
\end{description}
\end{center}
D'apr\`es la caract\'erisation des courbes lisses de $\text{Diff}(M)$, on en d\'eduit que $\widetilde{F}_{t}$ est une courbe lisse de $\text{Diff}(M)$ joignant $Id_{M}$ et $\widetilde{F}_{1}$ ce qui implique en particulier que $\widetilde{F}_{1}\in \text{Diff}^{\,0}(M)\,.$ De plus, $\vartheta(\widetilde{F}_{1},\Sigma_{0})=\widetilde{F}_{1}(\Sigma_{0})=F_{1}(\Sigma_{0})=\Sigma_{1}$ ce qui prouve le th\'eor\`eme.$\hfill\square$
\begin{remarque}  On pourrait montrer le Th\'eor\`eme \ref{theoreme 2} en utilisant le Th\'eor\`eme de Nash-Moser via le Th\'eor\`eme 2.4.1 de \cite{Hamilton}.
\end{remarque}
\begin{remarque} A partir du Th\'eor\`eme \ref{theoreme 2}\,, on peut montrer que la composante connexe $\big(Gr_{k}(M)\big)_{\Sigma}$ de la Grassmannienne est aussi homog\`eme sous l'action du groupe $\text{SDiff}(M,\mu)$ des diff\'eomorphimes de $M$ qui pr\'eservent une forme volume donn\'ee $\mu$ (voir \cite{Haller-Vizman})\,.
\end{remarque}

\section*{Appendice}
Dans cet appendice, on donne -- sans d\'emonstrations -- quelques r\'esultats techniques utiles pour la g\'eom\'etrie en dimension infinie, plus particuli\`erement pour l'\'etude des vari\'et\'es model\'ees sur des espaces de Fr\'echet (pour une introduction aux espaces de Fr\'echet, on pourra consulter \cite{Bierstedt-Bonet} ou \cite{Jarchow}, pour les vari\'et\'es model\'ees sur des espaces de Fr\'echet, \cite{Hamilton}, \cite{Kriegl-Michor}, etc.).
\begin{definition}
Soit F un espace de Fr\'echet, I un ouvert de $\mathbb{R}$ et $c\,:\,I\rightarrow F$ une application. On dit que c est d\'erivable sur I si pour tout $x\in I\,,$ le quotient $\big(c(x+h)-c(x)\big)/h$ converge lorsque $h\rightarrow 0\,,$ on note alors $c'$ sa d\'eriv\'ee. On dit qu'une courbe $c\,:\,I\rightarrow F$ est lisse si elle admet des d\'eriv\'ees \`a tous les ordres.
\end{definition}
La proposition ``folklorique'' suivante (voir \cite{Molitor2}) relie deux notions de calcul diff\'erentiel  sur les espaces de Fr\'echet. L'une utilise la notion de courbes lisses et est d\'evelopp\'ee dans \cite{Kriegl-Michor}, l'autre, plus classique, utilise la diff\'erentielle de G\^ateaux et est d\'evelopp\'ee, par exemple, dans \cite{Hamilton}, \cite{Milnor}, etc.
\begin{proposition}
Si $U\subseteq E$ est un ouvert d'un espace de Fr\'echet $E$ et $f\,:\,U\rightarrow F$ une application de $U$ dans un autre espace de Fr\'echet $F\,,$ alors $f$ est lisse (au sens de \cite{Hamilton}) si et seulement si $f\circ c$ est une courbe lisse de F pour toute courbe lisse $c\,:\,I\rightarrow U\,.$
\end{proposition}
Pour rendre cette derni\`ere proposition utile , nous avons besoin d'une bonne description (que l'on peut trouver dans \cite{Kriegl-Michor}) des courbes lisses de $\Gamma_{C^{\infty}}(M,\,E)$ o\`u $M$ est une vari\'et\'e compacte et $E\rightarrow M$ un fibr\'e vectoriel de rang fini (pour une description de la topologie de $\Gamma_{C^{\infty}}(M,\,E)\,,$ on pourra consulter \cite{dieudonne}, Proposition 17.2.2, page 238).
\begin{proposition}\label{caracterisation courbe lisse}
Si $s\,:\,I\rightarrow\Gamma_{C^{\infty}}(M,\,E)$ est une courbe lisse de $\Gamma_{C^{\infty}}(M,\,E)\,,$ alors l'application $s^{\wedge}\,:\, I\times M\rightarrow E,\,(t,x)\mapsto s_{t}(x)$ est une application lisse.\\
R\'eciproquement, si $f\,:\,I\times M\rightarrow E$ est une application lisse telle que $f(t,x)\in E_{x}$ pour tout $(t,x)\in I\times M\,,$ alors l'application $f^{\vee}\,:\,I\rightarrow \Gamma_{C^{\infty}}(M,\,E)$ d\'efinie par $f^{\vee}(t)(x):=f(t,x)$ est une courbe lisse de $\Gamma_{C^{\infty}}(M,\,E)\,.$
\end{proposition}
De cette proposition, on peut en d\'eduire facilement une caract\'erisation naturelle des courbes lisses des sous-vari\'et\'es de $C^{\infty}(M,N)$ pour laquelle on renvoie le lecteur \`a \cite{Kriegl-Michor}, Lemme 42.5, page 442.
\begin{definition}
Une suite $(x_{n})_{n\in\mathbb{N}}$ d'un espace de Fr\'echet $F$ ``converge rapidement'' vers $x\in F$ si pour tout $k\in \mathbb{N}\,,$ la suite $n^{k}(x_{n}-x)$ est born\'ee (born\'ee au sens des espaces topologiques localement convexes, voir \cite{Jarchow} ou \cite{Kriegl-Michor}).
\end{definition}
\begin{lemme}
Si $(x_{n})_{n\in\mathbb{N}}\,,$ est une suite d'un espace de Fr\'echet $F$ qui converge vers $x\in F\,,$ alors on peut trouver une sous-suite de $(x_{n})_{n\in\mathbb{N}}$ qui converge rapidement vers x.
\end{lemme}
De ce lemme ainsi que du ``special curve lemma'' de \cite{Kriegl-Michor}, page 16, on en d\'eduit
\begin{proposition}\label{special curve}
$\text{}$\,\,Si $(x_{n})_{n\in\mathbb{N}}$ est une suite d'un espace de Fr\'echet F qui converge vers $x\in F\,,$ alors (\`a sous-suite pr\`es) on peut trouver une courbe lisse $c\,:\,\mathbb{R}\rightarrow F$ telle que $c(\frac{1}{n})=x_{n}$ et $c(0)=x\,.$
\end{proposition}$\text{}$\\\\
$\text{}$\,\,\,\,\,\,\,\,\textbf{\begin{large}Remerciements.\end{large}}  Je tiens \`a remercier Tilmann Wurzbacher qui m'a encourag\'e \`a faire cet article et qui m'a chaleureusement accompagn\'e durant sa r\'edaction. 
%et ce, malgr\'e certains \'ev\`enements heureux.
\nocite{Ismagilov}
%\begin{footnotesize}\bibliography{biblio}\end{footnotesize}

\begin{thebibliography}{Ham82}

\bibitem[BB03]{Bierstedt-Bonet}
K.~D. Bierstedt and J.~Bonet.
\newblock \textit{Some aspects of the modern theory of {F}r\'echet spaces}.
\newblock {\em \emph{Rev. R. Acad. Cienc. Exactas F\'\i s. Nat. Ser. A Mat.}},
  97(2):159--188, 2003.

\bibitem[Die72]{dieudonne}
J.~Dieudonn{\'e}.
\newblock {\em \textit{Treatise on analysis. {V}ol. {III}}}.
\newblock {Academic Press}, New York, 1972.
\newblock Pure and Applied Mathematics, Vol. 10-III.

\bibitem[Ham82]{Hamilton}
R.~S. Hamilton.
\newblock \textit{The inverse function theorem of {N}ash and {M}oser}.
\newblock {\em \emph{Bull. Amer. Math. Soc. (N.S.)}}, 7(1):65--222, 1982.

\bibitem[Hir94]{Hirsch}
M.~W. Hirsch.
\newblock {\em \textit{Differential topology}}, volume~33 of {\em
  \emph{Graduate Texts in Mathematics}}.
\newblock {Springer-Verlag}, New York, 1994.

\bibitem[HV04]{Haller-Vizman}
S.~Haller and C.~Vizman.
\newblock \textit{Non-linear {G}rassmannians as coadjoint orbits}.
\newblock {\em \emph{Math. Ann.}}, 329(4):771--785, 2004.

\bibitem[Ism96]{Ismagilov}
R.~S. Ismagilov.
\newblock {\em \textit{Representations of infinite-dimensional groups}}, volume
  152 of {\em \emph{Translations of Mathematical Monographs}}.
\newblock {American Mathematical Society}, Providence, RI, 1996.

\bibitem[Jar81]{Jarchow}
H.~Jarchow.
\newblock {\em \textit{Locally convex spaces}}.
\newblock \emph{B. G. Teubner}, Stuttgart, 1981.

\bibitem[KM97]{Kriegl-Michor}
A.~Kriegl and P.~W. Michor.
\newblock {\em \textit{The convenient setting of global analysis}}, volume~53
  of {\em \emph{Mathematical Surveys and Monographs}}.
\newblock {American Mathematical Society}, Providence, RI, 1997.

\bibitem[Mil84]{Milnor}
J.~Milnor.
\newblock \textit{Remarks on infinite-dimensional {L}ie groups}.
\newblock In {\em Relativity, groups and topology, II (Les Houches, 1983)},
  pages 1007--1057. \emph{North-Holland}, Amsterdam, 1984.

\bibitem[Mol]{Molitor2}
M.~Molitor.
\newblock {\em \textit{Grassmanniennes non-lineaires, groupes de diff\'eomo-
  rphismes unimodulaires et quelques \'equations hamiltoniennes en dimension
  infinie}}.
\newblock \emph{Th\`ese de doctorat}, Universit\'e Paul Verlaine-Metz, France,
  2007.

\end{thebibliography}

\end{document}